\documentclass[preprint,10pt]{elsarticle}
\usepackage{amssymb}
\usepackage{amsthm,amsmath}

\journal{Linear Algebra and its Applications}

\newtheorem{thm}{Theorem}[section]
\newtheorem{prop}[thm]{Proposition}
\newtheorem{lem}[thm]{Lemma}
\newtheorem{cor}[thm]{Corollary}

\theoremstyle{definition}
\newtheorem{defn}[thm]{Definition}
\newtheorem{ex}[thm]{Example}

\theoremstyle{remark}
\newtheorem{remark}[thm]{Remark}

\newcommand{\R}{\mathbb{R}}  
\newcommand{\N}{\mathbb{N}}  
\newcommand{\Z}{\mathbb{Z}}  
\newcommand{\C}{\mathbb{C}} 
\newcommand{\s}{\{s_k\}_{k\in \N_0}} 
\newcommand{\ts}{\{t_k\}_{k\in \N_0}} 
\numberwithin{equation}{section}

\begin{document}
\begin{frontmatter}

\author{Hayoung Choi\corref{cor}}
\ead{hchoi2@uwyo.edu}

\author{Farhad Jafari\corref{}}
\ead{fjafari@uwyo.edu}

\cortext[cor]{Corresponding author}

\address{Department of Mathematics, University of Wyoming, Laramie, WY 82071, USA}

\title{Positive definite Hankel matrix completions and Hamburger moment completions}

\begin{abstract}
In this paper we give solutions to Hamburger moment problems with missing entries. The problem of completing partial positive sequences is considered. The main result is a characterization of positive definite completable patterns, namely patterns that guarantee the existence of Hamburger moment completion of a partial positive definite sequence. Moreover, several patterns which are not positive definite completable are given.
\end{abstract}

\end{frontmatter}


\section{Introduction}
A \emph{Hamburger moment sequence} of a given nondecreasing function $\sigma$ on $(-\infty,\infty)$ is defined by
\begin{equation*}
s_{k}:=\int_{-\infty}^{\infty} x^k \mathrm{d}\sigma(x) \quad \text{for all } k\in \N_{0},
\end{equation*}
assuming the integrals converge.

\begin{defn}
An infinite sequence $\s$ is \emph{positive} if the quadratic forms
\begin{equation*}
\sum_{i,j=0}^n x_i x_j s_{i+j} \geq 0 \quad \forall n\in
\N_0,~ x_0,\ldots,x_n \in \R.
\end{equation*}
\end{defn}

This is equivalent to the fact that Hankel matrices
\begin{equation}
H_n:=
\begin{bmatrix}
s_{0} & s_{1} & \cdots & s_{n} \\
s_{1} & s_{2} & \cdots & s_{n+1} \\
\vdots & \vdots & \ddots & \vdots \\
s_{n} & s_{n+1} & \cdots & s_{2n} \\
\end{bmatrix}
\end{equation}
are positive semidefinite for all $n \in \N_{0}$.

In 1920-21, Hamburger proved the the following remarkable theorem.
\begin{thm}[The classic Hamburger moment problem]
\label{intro:classic}
Given an infinite sequence of real numbers $\s$,
a necessary and sufficient condition for the existence of a nondecreasing $\sigma$ on $(-\infty,\infty)$ such that
\begin{equation*}
s_{k}=\int_{-\infty}^{\infty} x^k \mathrm{d}\sigma(x) \quad \text{for all } k\in \N_{0}
\end{equation*}
is that $\s$ is positive.
\begin{proof}
See (\cite{book:Akhiezer}, Theorem 2.1.1) or (\cite{book:completions}, Theorem 2.7.1).
\end{proof}
\end{thm}

\begin{defn}\label{defn:definite}
Let $\s$ be a positive sequence.
\begin{enumerate}
\item[i)] A sequence $\s$ is \emph{positive definite} if the Hankel matrix $H_n$ is positive definite for all $n\in \N_0$.\\
\item[ii)] A sequence $\s$ is \emph{positive semidefinite} if it is not positive definite.
\end{enumerate}
\end{defn}
Observe that if a sequence is positive semidefinite, then the Hankel matrix $H_n$ is positive semidefinite for all $n \in \N_0$ and at least one of them is singular. In fact, once one of the finite Hankel matrix is singular, all the following ones are also singular. Based on this definition, the sets of positive definite sequence and positive semidefinite sequence are mutually disjoint.\\

Numerous applications of moment problems, relying on a complete set of moments or truncated moment sequences, have been identified \cite{Las12}. For practical implementation of moment problems, it is vital to be able to deal with missing moment data since data obtained from physical experiments and phenomena are often corrupt or incomplete. Moment problems with missing entries are closely related to Hankel matrix completion problems (see \cite {book:completions} and \cite{Dominique}).

A \emph{partial sequence} is a sequence in which some terms are specified, while the remaining terms are unspecified and may be treated as free real variables.
A \emph{partial positive (semi)definite sequence} is a partial sequence if each of its fully specified principal submatrices of the Hankel matrix $H_n$ is positive (semi)definite for all $n \in \N_{0}$.
A \emph{Hamburger moment completion} (or a \emph{positive completion}) of a partial sequence is a specific choice of values for the unspecified terms resulting in a positive (semi)definite sequence. Note that a partial positive definite sequence may have either a positive definite completion or a positive semidefinite completion.
A \emph{Hamburger moment completion problem} asks whether a given partial sequence has a Hamburger moment completion.
A \emph{pattern} of a partial sequence is the set of positions of the specified entries. Denote the pattern of a partial sequence by the set of positive
integers
\begin{equation*}
P=\{k \in \N_0 ~:~ s_k \text{ is specified}\}.
\end{equation*}
We say that a pattern $P$ is \emph{positive (semi)definite completable} if every partial positive (semi)definite sequence with pattern $P$ has a positive completion, respectively.

Since every principal submatrix of a positive semidefinite matrix is positive semidefinite, it is trivial that a partial sequence has a positive completion only if it is a partial positive sequence.

\begin{ex}
The following is not a partial positive definite sequence.
\begin{equation}\label{intro:ex1}
 1,~ ?,~ \frac{1}{3},~ \frac{1}{7},~ \frac{1}{10},~?,~ ?,~ ?,~ ?,~ ,\cdots,
\end{equation}
where all entries only except $s_0$, $s_2$, $s_3$, and $s_4$ are unspecified, denoted by ?.
Since the principal submatrix of the Hankel matrix $H_2$ corresponding to \eqref{intro:ex1}
\begin{equation*}
\begin{bmatrix}
1 & \frac{1}{3}\\
\frac{1}{3} & \frac{1}{10}\\
\end{bmatrix}
\end{equation*}
has a negative determinant, the sequence \eqref{intro:ex1} has no positive definite completion.
\end{ex}

Thus, partial positivity is a necessary condition for the existence of a positive completion.
However, there exist a partial positive (semi)definite sequence which have no positive (semi)definite completion, respectively. In other words, partial positivity does not guarantee that there is a positive definite completion.
\begin{ex}
The following partial positive semidefinite sequence has no positive semidefinite completion:
\begin{equation}\label{intro:ex2}
 1,~ 1,~ 1,~ ?,~ 1,~ -1,~ ?,~ ?,~ ?,~ ?,\ldots,
\end{equation}
where only $s_0$, $s_1$, $s_2$, $s_4$, and $s_5$ are specified entries and denote missing entries by ?.
Let's check the Hankel matrix
\begin{equation*}
H_6=
\begin{bmatrix}
1 & 1 & 1 & ? & 1 &-1 & ?\\
1 & 1 & ? & 1 &-1 & ? & ?\\
1 & ? & 1 &-1 & ? & ? & ?\\
? & 1 &-1 & ? & ? & ? & ?\\
1 &-1 & ? & ? & ? & ? & ?\\
-1& ? & ? & ? & ? & ? & ?\\
? & ? & ? & ? & ? & ? & ?\\
\end{bmatrix}.
\end{equation*}
Since fully specified principal submatrices of Hankel matrices corresponding to \eqref{intro:ex2} are only $H_6[\{1,2\}]$ and $H_6[\{1,3\}]$, which are positive semidefinite, the partial sequence \eqref{intro:ex2} is a partial positive semidefinite. However, since principal submatrices $H_6[\{1,2,3\}]$ and $H_6[\{2,3,4\}]$ cannot be positive semidefinite simultaneously, it has no positive semidefinite completion. Thus, the pattern $P=\{0,1,2,4,5\}$ is not positive semidefinite completable.
\end{ex}

\begin{ex}
Note that the partial sequence
\begin{equation}\label{intro:ex3}
 1,~ ?,~\frac{1}{2},~ ?, ~\frac{1}{3},~?, ~\frac{1}{4},~?, ~\frac{1}{5},~?, ~\frac{1}{6},~?, \cdots.
\end{equation}
is partial positive definite since the Hilbert matrix is totally positive (all of its minors are positive).
Since the subsequence $\{\frac{1}{k+1}\}_{k=0}^{\infty}$ is positive definite, by Theorem \ref{intro:classic}
there exists a nondecreasing function $\sigma(x)$ on $(-\infty, \infty)$ such that
\begin{equation}\label{ex1:submoment}
\frac{1}{k+1}=\int_{-\infty}^{\infty} x^k \mathrm{d}\sigma(x) \quad \textup{for all } k\in \N_0.
\end{equation}
However, the sequence \eqref{ex1:submoment} is not a positive definite completion of \eqref{intro:ex3}, i.e., it does not mean that
there exists a nondecreasing function $\widetilde{\sigma}(x)$ on $(-\infty, \infty)$ such that
\begin{equation*}
\frac{1}{k+1}=\int_{-\infty}^{\infty} x^{2k} \mathrm{d}\tilde{\sigma}(x) \quad
\textup{for all } k\in \N_0.
\end{equation*}
In fact, by Corollary \ref{cor:mainresult2} the partial sequence \eqref{intro:ex3} has a positive definite completion.
\end{ex}

In this paper we characterize the partial positive (semi)definite sequences that admit positive (semi)definite completions.
After a series of preliminary definitions and results in Section 2 that will be used throughout this paper, in Section 3, we identify patterns that have positive definite Hankel matrix completions. The main result of this section shows that a pattern that is a subset of odd positive integers is positive definite completable.
These results suggest that positive definite completable sequences can be quite sparse.
Also, any fully specified truncated pattern is positive definite completable.
In Section 4,
we show that any pattern that is an arithmetic progression, with an even offset, is positive semidefinite completable. Surprisingly, to prove this result our key observation reduced the Hamburger moment completion problem to the Stieltjes moment completion problem, and using an integral representation, we were able to construct the completion of the sequence.

In Section 5, we show that if a sequence is positive semidefinite completable, then it is positive definite completable, but not conversely.
Also, we characterize partial positive sequences that have positive definite completions but not positive semidefinite completions.  Finally, by considering $3 \times 3 $ and $4 \times 4$ positive definite completable Hankel matrices, in Section 6, we characterize general partial Hankel matrices that do not have positive definite Hankel matrix completions.

\section{Preliminaries}
In this section, we gather a series of results that will be used throughout this paper. While primarily focused on the Hamburger moment problem,  the Stieltjes and the trigonometric (Carath\'{e}odory) moment problems make appearances in the analysis. We recall that if the function $ \sigma (x) $ is supported on the nonnegative real line and on $ [-\pi, \pi] $, the moment problem is the Stieltjes moment problem and the Carath\'{e}odory moment problem, respectively. Two fundamental theorems of Schur provide criteria for positivity that will be used repeatedly throughout the results sections. It is worthwhile to point out that there is a significant difference between the definite and semidefinite character of moment sequences in their integral representations.
\begin{thm}[\cite{book:widder}, Theorem 12a]\label{moment}
A necessary and sufficient condition that there exists a non-decreasing function $\sigma(x)$ with infinitely many points of increase (with a finite number of points of increase) such that
\begin{equation}
s_{k}=\int_{-\infty}^{\infty} x^k \mathrm{d}\sigma(x) \quad \textup{for all } k\in \N_{0}
\end{equation}
is that the sequence $\s$ is positive definite (semidefinite), respectively.
\end{thm}

The following statement is a special case of the Hamburger problem which is called the Stieltjes moment problem.
\begin{thm}[The classic Stieltjes moment problem]\label{Prelim:stieltjes}
A necessary and sufficient condition that there exists a non-decreasing function $\sigma(x)$ with infinitely many points of increase (with a finite number of points of increase) such that
\begin{equation}
s_{k}=\int_{0}^{\infty} x^k \mathrm{d}\sigma(x) \quad \textup{for all } k\in \N_{0}
\end{equation}
is that the sequences $\s$ and $\{s_k\}_{k\in\N}$ are positive definite (the sequences $\s$ or $\{s_k\}_{k\in\N}$  is all positive, at least one of them being positive semidefinite), respectively.
\begin{proof}
See (\cite{book:widder}, Theorem 13a, 13b, and 13c)
\end{proof}
\end{thm}
Note that in addition to the positive definiteness of the sequence $\s $, the auxiliary condition of positive definiteness of the shifted sequence $\{s_k\}_{k\in\N}$ is added here.
\newline

Let $\mathbf{PSD}_n$ be the set of all $n\times n$ positive semidefinite matrices. If $A,B\in \mathbf{PSD}_n$, then $A\circ B \in \mathbf{PSD}_n$, where
$\circ$ is the Hadamard product or the Schur product (See \cite{book:Matrixanalysis}, Theorem 7.5.3).
Furthermore, any nonnegative linear combination of positive semidefinite matrices is also positive semidefinite.
Thus, $\mathbf{PSD}_n$ is a positive cone in the vector space of all Hermitian matrices (See \cite{book:completions}, Lemma 1.1.3). Therefore, the following statements are trivial.

\begin{prop}\label{positive:thm1}
$(1)$ If $\{s_n\}$ and $\{t_n\}$ are positive (semi)definite sequences, then $\{\alpha s_n + \beta t_n\}$ is a positive (semi)definite sequence for $\alpha>0$, $\beta>0$ respectively;\\
$(2)$ Term-wise product of two positive (semi)definite sequences is a positive (semi)definite sequence, respectively.\\
\end{prop}
It follows that if $\s$ is positive (semi)definite and if
$f(x)=\sum_{n=0}^\infty c_n x^n$ with $c_n\geq 0$ for all $n$, then
so is $\{f(s_k)\}_{k=0}^{\infty}$.
The converse of this theorem is due to Rudin (\cite{absmonoto}).

Note that $H_n$ is positive semidefinite if and only if it is Hermitian and all its principal minors are nonnegative. Given a Hankel matrix $H_n$ corresponding to a positive sequence $\s$, consider a principal submatrix
$H_n[\alpha]$, where $\alpha=\{\alpha_0,\alpha_1,\alpha_2,\cdots,\alpha_r\}\subseteq\N_0$. Then, $H_n[\alpha]$ is positive semidefinite. However, $H_n[\alpha]$ may not be a Hankel matrix. For $H_n[\alpha]$ to be a Hankel matrix we require that for any $k$,
\begin{equation}
\alpha_{k}=\frac{\alpha_{k-1}+\alpha_{k+1}}{2}
\end{equation}
and
$H_n[\alpha]$ has (1,1) entry $s_{2\alpha_0}$. Thus, one gets the following theorem.

\begin{thm}\label{thm:subsequence}
Let $\s$ be a positive sequence. Then, a subsequence $\{s_{l_k}\}_{k\in \N_0}$ of $\s$
is positive if
\begin{equation}
l_k = kd + l_0, \label{eq:subsequence}
\end{equation}
where $d\in \N$ and $l_0 \in 2 \N_{0}~$.
\end{thm}
However, there are, of course, many positive subsequences that are not of the form as in Theorem \ref{thm:subsequence}.
For example, since $\{\frac{1}{k+1}\}_{k=0}^{\infty}$ is a positive definite sequence, by Proposition \ref{positive:thm1} (2) its subsequence $\{\frac{1}{(k+1)^2}\}_{k=0}^{\infty}$ is positive definite,
which is not of the form (\ref{eq:subsequence}).

\begin{defn}
A \emph{trigonometric moment sequence} of a nondecreasing function $\sigma(\theta)$ on $[-\pi,\pi]$ is defined by
\begin{equation*}
t_{k}=\frac{1}{2\pi}\int_{-\pi}^{\pi} e^{ik\theta} \mathrm{d}\sigma(\theta) \quad \text{for all } k\in \Z,
\end{equation*}
assuming the integrals converge.
\end{defn}

The classical trigonometric moment problem states that given an infinite sequence of complex numbers $\{t_k\}_{k\in \Z}$, a necessary and sufficient condition for the existence of a nondecreasing function $\sigma(\theta)$ on $[-\pi,\pi]$ such that
\begin{equation*}
t_{k}=\frac{1}{2\pi}\int_{-\pi}^{\pi} e^{ik\theta} \mathrm{d}\sigma(\theta) \quad \text{for all } k\in \Z,
\end{equation*}
is that all Toeplitz matrices
\begin{equation*}
T_n:=\begin{bmatrix}
t_{0} & \bar{t}_{1} & \cdots & \cdots & \bar{t}_{n} \\
t_{1} & t_{0} & \bar{t}_{1} & \cdots & \bar{t}_{n-1} \\
\vdots & \ddots & \ddots & \ddots & \vdots\\
\vdots & \ddots & \ddots & \ddots & \bar{t}_{1}\\
t_{n} & t_{n-1} & \cdots & t_{1} & t_{0} \\
\end{bmatrix}
\end{equation*}
are positive semidefinite for all $n \in \N_{0}$.


In 1981, H. Dym and I. Gohberg studied extensions of band matrices with band inverses \cite{DG81}.
In 1984, R. Grone, C. R. Johnson, E. M. S\'{a}, and H. Wolkowicz considered positive definite completion of partial Hermitian matrices. They showed that if the undirected graph of the specified entries is chordal, a positive definite completion necessarily exists \cite{Grone}. C. Johnson studied positive definite Toeplitz matrix completions in 1997. In \cite{Charles} he proved that a pattern $P$ of an $(n+1)\times(n+1)$ partial Toeplitz matrix is positive (semi)definite completable if and only if $P=\{k,2k,\ldots,mk\}$
for some $m\in \N$ and $k\in \N$.
M. Bakonyi and G. Naevdal gave a characterization of certain class of subsets that guarantee the extension property for trigonometric moments by using matrix extension methods \cite{Bakonyi}.
Also, they characterized all finite subsets of $\Z^2$ which is positive (semi)definite completable \cite{BN00}. Moreover, Bakonyi, Rodman, Spitkovsky, and Woerdeman proved that each infinite band of a certain form is positive (semi)definite completable \cite{BRSW96}.
\newline

The following two facts will be used for the main results of this paper.
\begin{prop}\label{prop:schur}
Let $M=
\begin{bmatrix}
A & B \\
B^T & C \\
\end{bmatrix}$
be a symmetric matrix such that $A$ is invertible and $S$ is the Schur complement of $A$ in $M$, which is
$S=C-B^TA^{-1}B$.
Then,
\begin{equation*}
M>0~ \Leftrightarrow~ A>0 \textup{ and } S>0.
\end{equation*}
\end{prop}

\begin{thm}[Schur complements and determinant formula]\label{thm:schurcomp}
Let $A=[a_{ij}]\in M_n(\C)$ be given and $\alpha\subseteq \{1,\cdots,n\}$ be
an index set such that $A[\alpha]$ is nonsingular. Then,
\begin{equation*}
\det{A}=\det{A[\alpha]}\det{\big(A[\alpha^c]-A[\alpha^c,\alpha]A[\alpha]^{-1}A[\alpha,\alpha^c]\big)}.
\end{equation*}
\end{thm}

\section{Positive Definite Hankel Completions}

A \emph{partial Hankel matrix} is a partial matrix that is a Hankel matrix to the extent to which it is specified; that is, all specified entries lie along certain specified skew-diagonals (the diagonals perpendicular to the main diagonals).
A \emph{positive (semi)definite completion} of a partial Hankel matrix is a specific choice of values for unspecified entries resulting in a positive semidefinite Hankel matrix.
A partial matrix is \emph{partial positive (semi)definite} if each of its fully specified principal submatrices is positive (semi)definite.
A \emph{pattern} of a partial Hankel matrix is the set of positions of specified entries.
Denote the pattern of an $(n+1)\times(n+1)$ partial Hankel matrix $H_n$ by
\begin{equation*}
P=\{k \in \{0,1,2,\dots,2n\} ~:~ s_k \text{ is specified}\}.
\end{equation*}
We say that a pattern $P$ of an $(n+1)\times(n+1)$ partial Hankel matrix $H_n$ is \emph{positive (semi)definite  completable} if every partial positive (semi)definite Hankel matrix defined on the pattern $P$ has a positive completion.

In the next series of results, we show that if a pattern $P$ is any subset of odd integers, then the pattern $P$ is positive definite completable. In particular, if $P$ is any subset of all prime numbers, then the pattern $P$ is positive definite completable.
Also, we show that while $P = 2 \N_0$ is positive definite completable, its subsets may not be completable.

If $H_n=[s_{i+j}]_{i,j=0,\ldots,n}$ is an $(n+1)\times(n+1)$ (partial) Hankel matrix with entries $s_0,s_1,\ldots,s_{2n}$, we use $H_n$, $H(s_0,s_1,\ldots,s_{2n})$, or $H$ interchangeably.

\begin{lem}\label{lem:oddgiven}
Let $H(s_0,s_1,\ldots,s_{2n+1},?)$ be a partial positive definite Hankel matrix with the only missing entry $s_{2n+2}$, denoted by $?$. Then,
it has a positive definite Hankel completion.
\begin{proof}
Partition $H_{n+1}:=H(s_0,s_1,\ldots,s_{2n+1},?)$
as
\begin{equation*}
H_{n+1}=
\begin{bmatrix}
H_n & v_n \\
v_n^T & s_{2n+2}
\end{bmatrix},
\end{equation*}
where
$v_n^T=
\begin{bmatrix}
s_{n+1} & s_{n+2} & \cdots & s_{2n+1}
\end{bmatrix}$.
Since $H_{n+1}$ is a partial positive definite,
every fully specified principal submatrix is positive definite matrix, so $H_n>0$.
By Proposition \ref{prop:schur} the following holds:
\begin{equation*}
H_{n+1}>0 ~ \Longleftrightarrow~
s_{2n+2}-v_n^T H_n^{-1} v_n >0.
\end{equation*}
Choosing $s_{2n+2}>0$ such that $s_{2n+2}>v_n^T H_n^{-1} v_n$,
it follows that $H_{n+1}>0$.
\end{proof}
\end{lem}

\begin{lem}\label{lem:pricipal}
Let $H_{n+1}:=H(s_0,s_1,\ldots,s_{2n},?,s_{2n+2})$ be a partial Hankel matrix with the only missing entry $s_{2n+1}$, denoted by $?$.
A necessary and sufficient condition that $H_{n+1}$ is a partial positive definite matrix is
\begin{equation*}
H_n>0 \text{  and }
s_{2n+2}-w_n^T H_{n-1}^{-1} w_n >0,
\end{equation*}
where $w_n^T=
\begin{bmatrix}
s_{n+1} & s_{n+2} & \cdots & s_{2n} \\
\end{bmatrix}$.
\begin{proof}
Consider the matrix
\begin{equation*}
H_{n+1}=
\begin{bmatrix}
s_{0} & s_{1} & \cdots & s_{n-1} & s_{n} & s_{n+1} \\
s_{1} & s_{2} & \cdots & s_{n} & s_{n+1} & s_{n+2} \\
\vdots & \vdots & \ddots & \vdots & \vdots & \vdots\\
s_{n-1} & s_{n} & \cdots & s_{2n-2} & s_{2n-1} & s_{2n} \\
s_{n} & s_{n+1} & \cdots & s_{2n-1} & s_{2n} & ? \\
s_{n+1} & s_{n+2} & \cdots & s_{2n} & ? & s_{2n+2} \\
\end{bmatrix}.
\end{equation*}
Note that all principal submatrices of $H_{n+1}$
can be classified with respect to
the last two columns and rows as follows:
\begin{equation*}
H_{n+1}[\alpha],~H_{n+1}[\alpha \cup \{n+1\}],~
H_{n+1}[\alpha \cup \{n+2\}],~H_{n+1}[\alpha \cup \{n+1,n+2\}]
\end{equation*}
for any $\alpha \subseteq \{1,\cdots,n\}$.
Note that both of $H_{n+1}[\alpha]$ and $H_{n+1}[\alpha \cup \{n+1\}]$ are principal submatrices of $H_n$ and that
$H_{n+1}[\alpha \cup \{n+1,n+2\}]$ is not a fully specified principal submatrix.
Thus, all principal submatrices of $H_{n+1}$ are positive definite
if and only if $H_{n}>0$ and $H_n[\{1,\cdots,n,n+2\}]>0$.
By Proposition \ref{prop:schur} it follows that
\begin{equation*}
H_n[\{1,\cdots,n,n+2\}]>0~
\Longleftrightarrow~
s_{2n+2}-w_n^T H_{n-1}^{-1} w_n >0.
\end{equation*}
\end{proof}
\end{lem}

\begin{lem}\label{lem:evengiven}
Let $H(s_0,\ldots,s_{2n},?,s_{2n+2})$ be a partial positive definite Hankel matrix with the only missing entry $s_{2n+1}$, denoted by $?$. Then, it has a positive definite Hankel completion.
\begin{proof}
Partition $H_{n+1}:=H(s_0,s_1,\ldots,s_{2n},?,s_{2n+2})$ as
\begin{equation*}
H_{n+1}=
\begin{bmatrix}
H_{n-1} & v_{n-1} & w_n\\
v_{n-1}^T & s_{2n} & s_{2n+1}\\
w_{n}^T & s_{2n+1} & s_{2n+2}\\
\end{bmatrix}.
\end{equation*}
By Proposition \ref{prop:schur}, it follows that
\begin{align*}
H_{n+1}>0
&\Longleftrightarrow
\begin{bmatrix}
s_{2n} & s_{2n+1}\\
s_{2n+1} & s_{2n+2}\\
\end{bmatrix}
-
\begin{bmatrix}
v_{n-1}^T \\
w_{n}^T\\
\end{bmatrix}
H_{n-1}^{-1}
\begin{bmatrix}
v_{n-1} & w_{n}
\end{bmatrix}>0\\
&\Longleftrightarrow
\begin{bmatrix}
s_{2n} - v_{n-1}^T H_{n-1}^{-1}v_{n-1} & s_{2n+1} - v_{n-1}^T H_{n-1}^{-1}w_{n}\\
s_{2n+1} - w_{n}^T H_{n-1}^{-1}v_{n-1} & s_{2n+2} - w_{n}^T H_{n-1}^{-1}w_{n}\\
\end{bmatrix}
>0
\end{align*}
Since $H_n>0$, $s_{2n} - v_{n-1}^T H_{n-1}^{-1}v_{n-1}>0$. By Lemma \ref{lem:pricipal},
$s_{2n+2} - w_{n}^T H_{n-1}^{-1}w_{n}>0$. Take $s_{2n+1}= v_{n-1}^T H_{n-1}^{-1}w_{n}$.
\end{proof}
\end{lem}

\begin{thm}\label{thm:oddevenmissing}
Let $H(s_0,\ldots,s_{2n},?,?)$ be a partial positive definite Hankel matrix with the only missing entries $s_{2n+1}$ and $s_{2n+2}$, denoted by $?$. Then, it has a positive definite Hankel completion.
\begin{proof}
Since $H(s_0,\ldots,s_{2n},?,?)$ is partial positive definite, $H_{n-1}>0$, implying $w_{n}^T H_{n-1}^{-1}w_{n}>0$. Then one can find $s_{2n+2}>0$ such that $s_{2n+2} - w_{n}^T H_{n-1}^{-1}w_{n}>0$. Then, by Lemma \ref{lem:pricipal} the partial Hankel matrix $H(s_0,\ldots,s_{2n},?,s_{2n+2})$ is partial positive definite. Thus, by Lemma \ref{lem:evengiven} there is $s_{2n+1}\in \R$ such that $H(s_0,\ldots,s_{2n},s_{2n+1},s_{2n+2})>0$.
\end{proof}
\end{thm}

\begin{thm}\label{thm:oddsubset}
A pattern $P\subset 2\N_0+1$ is positive definite completable.
\begin{proof}
Let $P\subset 2\N_0+1$ be any pattern and $\s$ be a partial positive sequence with the pattern $P$.
Note that any Hankel matrix $H_n$ defined on $P$ has no fully specified principal submatrices for $n\in \N_0$. Assume $s_0=1$. When $1\in P$, i.e., $s_1$ is specified, then one can find $s_2>0$ such that $H_1>0$. When $1\notin P$, choose any arbitrary $s_1\in \R$ and then find $s_2>0$ such that $H_1>0$.
Suppose that $s_0$, $s_1$, $s_2$, $\cdots$, $s_{2n}$ are given such that $H_n(s_0,s_1,\ldots,s_{2n})>0$ for some $n\in\N$. There are two different possibilities:
(i) When $2n+1\in P$, since $H_{n+1}(s_0,s_1,\ldots,s_{2n},s_{2n+1},?)$ is a partial positive definite Hankel matrix, by Lemma \ref{lem:oddgiven} there exists $s_{2n+2}>0$ such that $H_{n+1}>0$.
(ii) When $2n+1\notin P$,
Since $H_{n+1}(s_0,s_1,\ldots,s_{2n},?,?)$ is a partial positive Hankel matrix, by Theorem \ref{thm:oddevenmissing} there exists $s_{2n+1}\in \R$ and $s_{2n+2}>0$ such that $H_{n+1}>0$.
Thus, there exists $s_0$, $s_1$, $\cdots$, $s_{2n+2}$ such that $H_{n+1}(s_0,s_1,\ldots,s_{2n},s_{2n+1},s_{2n+2})>0$.
Thus, by induction there exists a positive definite sequence $\s$. Therefore, the pattern $P$ is positive definite completable.
\end{proof}
\end{thm}

It is easy to extend Theorem \ref{thm:oddsubset} somewhat.
\begin{cor}\label{cor:union}
Let $P\subset 2\N_0+1$. Then the pattern $P\cup \{0,1,\ldots,m\}$
is positive definite completable.
\begin{proof}
Let $\s$ be a partial positive sequence with the pattern $P\cup \{0,\ldots,m\}$.
(i) When $m=2n+1$, the corresponding Hankel matrix $H_{n+1}(s_0,s_1,\ldots,s_{2n+1},?)$ is partial positive definite. By Lemma \ref{lem:oddgiven} it follows that there exists $s_{2n+2}>0$ such that $H_{n+1}>0$. Next there are two different possibilities: $2n+3\in P$
and $2n+3\notin P$. When $2n+3\in P$, since $H_{n+2}(s_0,s_1,\ldots,s_{2n+2},s_{2n+3},?)$ is a partial positive definite Hankel matrix, by Lemma \ref{lem:oddgiven} there exists $s_{2n+4}>0$ such that $H_{n+2}>0$. When $2n+3\notin P$,
since $H_{n+2}(s_0,s_1,\ldots,s_{2n+2},?,?)$ is a partial positive Hankel matrix, by Theorem \ref{thm:oddevenmissing} there exists $s_{2n+3}\in \R$ and $s_{2n+4}>0$ such that $H_{n+2}>0$.
Thus, there exists $s_0$, $s_1$, $\cdots$, $s_{2n+4}$ such that $H_{n+2}(s_0,s_1,\ldots,s_{2n+3},s_{2n+4})>0$.
Thus, by induction there exists a positive definite completion.

(ii) When $m=2n$, if $2n+1\in P$ then by Lemma \ref{lem:oddgiven} there exists $s_{2n+2}>0$ such that $H_{n+1}>0$. If $2n+1\notin P$
then by Theorem \ref{thm:oddevenmissing} there exists $s_{2n+1}\in \R$ and $s_{2n+2}>0$ such that $H_{n+1}>0$. By induction there exists a positive definite completion.

Therefore, by (i) and (ii) the pattern $P\cup \{0,1,\ldots,m\}$ is positive definite completable.
\end{proof}
\end{cor}

The following example is quite interesting and suggests that patterns can be quite complicated.
\begin{cor}\label{cor:prime}
Let $P$ be the set of all prime numbers. Then the pattern $P$ is positive definite completable.
\begin{proof}
Note that $P-\{2\}\subset 2\N_0+1$. Let $\s$ be a partial positive sequence with the pattern $P$.
Assume $s_0=1$. Since $s_2>0$, there exists $s_1\in \R$ such that $H_1(s_0,?,s_2)>0$. Since $H_2(s_0,s_1,s_2,s_3,?)$ is a partial positive definite Hankel matrix, by Lemma \ref{lem:oddgiven} one can choose $s_4>0$ such that $H_2>0$. Apply the proof of Theorem \ref{thm:oddsubset}.
\end{proof}
\end{cor}

\begin{remark}
By Theorems \ref{thm:oddsubset} and \ref{thm:mainresult} the pattern $P=2\N_0+1$ and $P= 2\N_0$ are positive definite completable. Moreover,
the patterns $\{2,4,\ldots,2n\}$ and
$\{1,3,\ldots,2n+1\}$ are positive definite completable as well.
However, a pattern $P\subsetneq 2\N_0$ may or may not be positive definite completable.
For example, the pattern $P=\{0,2,8\}\subset 2\N_0$ is not positive definite completable. By Theorem \ref{thm:3x3main} the following partial positive definite sequence has no positive definite completion.
\begin{equation*}
1,~?,~\frac{1}{2},~?,~?,~?,~?,~?,~\frac{1}{16}, ~?,~?,~?,~\dots,
\end{equation*}
where only $s_0$, $s_2$, and $s_8$ are specified entries and ? denotes missing entries.
\end{remark}

\begin{cor}\label{cor:trunc}
The pattern $P=\{0,1,2,\cdots,m\}$ is positive definite completable.
\begin{proof}
Case 1: When $m=2n+1$, let
$\s$ be a partial positive definite sequence with the pattern $P$.
Since $H_{n+1}(s_0,s_1,\ldots,s_{2n+1},?)$ is partial positive definite, by
Lemma \ref{lem:oddgiven} $H_{n+1}$ has a positive definite completion, say $s_{2n+2}$.
Then the partial Hankel matrix $H_{n+2}(s_0,s_1,\ldots,s_{2n+2},?,?)$ is partial positive definite. By Theorem \ref{thm:oddevenmissing} there exist $s_{2n+3}\in \R$ and $s_{2n+4}>0$ such that $H_{n+2}>0$.
Continue this process by induction.

Case 2: When $m=2n$,
let $\s$ be a partial positive definite sequence with the pattern $P$.
Since $H_{n+1}(s_0,s_1,\ldots,s_{2n},?,?)$ is partial positive definite, by Theorem \ref{thm:oddevenmissing} there exist $s_{2n+1}\in \R$ and $s_{2n+2}>0$ such that $H_{n+1}>0$. Continuous this process.
\end{proof}
\end{cor}

\begin{remark}\label{remark:semicompletion}
By Corollary \ref{cor:trunc} for the given real numbers $s_0$, $s_1$, $\ldots$, $s_m$,
a sufficient condition for the existence of a nondecreasing function $\sigma$ on $(-\infty,\infty)$ such that
\begin{equation*}
s_{k}=\int_{-\infty}^{\infty} x^k \mathrm{d}\sigma(x) \quad \text{for all } 0\leq k \leq m
\end{equation*}
is that the partial sequence $s_0$, $s_1$, $\ldots$, $s_m$, $?$, $?$, $\ldots$ is partial positive definite.
However, the $P=\{0,1,2,\cdots,m\}$ is not positive semidefinite completable. M. Bakonyi and H. Woerdeman give an extra condition for the pattern $P$ to have a positive semidefinite completion (\cite{book:completions}, Theorem 2.7.5).
\end{remark}

\section{Hamburger Moment Completions}

For positive semidefinite completable patterns, arithmetic progression patterns play a crucial role. The following theorem is the main result of this section.
\begin{thm}\label{thm:mainresult}
If the pattern $P = d\N_0+l_0$ for some $d\in \N$ and $l_0 \in 2\N_{0}$, then $P$ is positive semidefinite completable.
\begin{proof}
Let $P=d\N_0+l_0$ be the pattern with $d\in \N$ and $l_0 \in 2\N_{0}$.
Suppose that $\s$ is a partial positive semidefinite sequence with the pattern $P$. Then, by Theorem \ref{thm:subsequence} its subsequence
$\{s_{dk+l_0}\}_{k\in\N_0}$
is a positive semidefinite sequence.
By Theorem \ref{moment}, there exists a non-decreasing function $\sigma(x)$ $(-\infty<x<\infty)$ with finitely many points of increase such that
\begin{equation}\label{subint}
s_{dk+l_0}=\int_{-\infty}^{\infty} x^k d\sigma(x) \qquad \textit{for all } k \in\N_0.
\end{equation}
Let $\lambda_0 < \lambda_1 < \ldots< \lambda_n $ be the points of increase of $\sigma(x)$ with
\begin{equation*}
\mu_i := \sigma(\lambda_i + 0) - \sigma(\lambda_i -0) \quad \textit{ for all }
0\leq i\leq n.
\end{equation*}
Without loss of generality we assume $\lambda_i\neq 0$.
If we exclude the point of increase $\lambda_i=0$ for some $0 \leq i \leq n$, then the function $\sigma$ may change, but the integrals \eqref{subint} remain invariant.
Then, by the definition of Stieltjes integral
\begin{equation}\label{eq:seq_def_stieltjes}
s_{dk+l_0}=\sum_{i=0}^{n} \mu_i \lambda_i^k \quad \textit{for all } k\in\N_0.
\end{equation}
Define $\phi :\R \rightarrow \R$ by 
\begin{equation*}
\phi(x) = \sqrt[d]{x}.
\end{equation*}
Define $\widetilde{\sigma}:\R \longrightarrow \R$ as a non-decreasing function with finitely many points of increase,
\begin{equation}\label{philambdas}
\phi(\lambda_0) < \phi(\lambda_1) < \phi(\lambda_2) < \cdots < \phi(\lambda_n)
\end{equation}
with
\begin{equation*}
\eta_i:= \frac{\mu_i} {\phi(\lambda_i)^{l_0}} 
= \widetilde{\sigma}(\phi(\lambda_i) + 0) - \widetilde{\sigma}(\phi(\lambda_i) -0) \quad \textit{for all } 0\leq i\leq n.
\end{equation*}
Since $l_0\in 2\N_0$, $\eta_i>0$ for all $0\leq i\leq n$.
When $d\in 2\N_0+1$, $\phi(\lambda_i)$ is well-defined.
Since $\phi$ is an increasing function, \eqref{philambdas} holds, i.e., the function $\widetilde{\sigma}$ is well-defined. When $d\in 2\N$, not only by hypothesis the sequence $\{s_{dk+l_0}\}_{k\in\N_0}$ is positive semidefinite, but also the sequence $\{s_{dk+l_0}\}_{k\in\N}$ is positive semidefinite. By Theorem \ref{Prelim:stieltjes} it follows that the points of increase of $\sigma(x)$ are all positive, i.e., $\lambda_i > 0$ for all $i$. Hence, $\phi(\lambda_i)$ is well-defined and \eqref{philambdas} holds.
Define the sequence $\{\tilde{s}_k\}_{k\in\N_0}$ by
\begin{equation}
\tilde{s}_k:=\int_{-\infty}^{\infty} x^k d\widetilde{\sigma}(x) \quad \textit{for all }
k\in\N_0.
\end{equation}
We claim that $\{\tilde{s}_k\}_{k\in\N_0}$ is a positive semidefinite completion of the partial sequence $\s$.
Indeed, by the definition of Stieltjes integral one gets that
\begin{equation}
\tilde{s}_k=\sum_{i=0}^{n} \eta_i (\phi(\lambda_i))^k.
\end{equation}
By (\ref{eq:seq_def_stieltjes}) it follows that $\tilde{s}_{dk+l_0}=s_{dk+l_0}$ for all $k\in \N_0$.
Hence, $\s$ has a positive semidefinite completion.
\end{proof}
\end{thm}

It is interesting to note that in the case when $ d $ was even in the above theorem, the Hamburger moment completion problem is converted to a Stieltjes moment completion problem. Here one readily gets that not only the original sequence is positive semidefinite, but the shifted sequence is also positive semidefinite.

The following corollary shows that if $l_0 = 0 $, then the partial positive pattern $P$ is actually positive definite completable.
The proof is very similar to the previous proof, but it is sufficiently different that for clarity we present it here.

\begin{cor}\label{cor:mainresult2}
If the pattern $P = d\N_0$ for some $d\in \N_0$, then $P$ is positive definite completable.
\begin{proof}
Let $P=d\N_0$ be the pattern with $d\in \N$.
Suppose that $\s$ is a partial positive definite sequence with the pattern $P$. Then, by Theorem \ref{thm:subsequence} its subsequence
$\{s_{dk}\}_{k\in\N_0}$
is a positive definite sequence.
By Theorem \ref{moment}, there exists a non-decreasing function $\sigma(x)$ $(-\infty<x<\infty)$ with infinitely many points of increase such that
\begin{equation}
s_{dk}=\int_{-\infty}^{\infty} x^k d\sigma(x) \quad \textit{for all } k \in\N_0.
\end{equation}
Let $\lambda_0 < \lambda_1 < \ldots< \lambda_n < \cdots$ be the points of increase of $\sigma(x)$ with
\begin{equation*}
\mu_i := \sigma(\lambda_i + 0) - \sigma(\lambda_i -0) \quad \textit{for all }
 i\in\N_0.
\end{equation*}
Without loss of generality we assume $\lambda_i\neq 0$.
Then, by the definition of Stieltjes integral
\begin{equation}\label{eq:seq_def_stieltjes2}
s_{dk}=\sum_{i=0}^{\infty} \mu_i \lambda_i^k \quad \textit{for all } k\in\N_0.
\end{equation}
Define $\phi :\R \rightarrow \R$ by 
\begin{equation*}
\phi(x) = \sqrt[d]{x}.
\end{equation*}
Define $\widetilde{\sigma}:\R \longrightarrow \R$ as a non-decreasing function with infinitely many points of increase at
\begin{equation}\label{philambdas2}
\phi(\lambda_0) < \phi(\lambda_1) < \phi(\lambda_2) < \cdots < \phi(\lambda_n) <\cdots
\end{equation}
with
\begin{equation*}
\mu_i 
= \widetilde{\sigma}(\phi(\lambda_i) + 0) - \widetilde{\sigma}(\phi(\lambda_i) -0) \quad \textit{for all } n\in \N_0.
\end{equation*}
Referring to the definition of $ \eta_i $ in Theorem \ref{thm:mainresult}, since $0=l_0\in 2\N_0$, $\eta_i > 0 $ for all $i\in \N_0$.
When $d\in 2\N_0+1$, $\phi(\lambda_i)$ is well-defined and
since $\phi$ is an increasing function, \eqref{philambdas2} holds, i.e., the function $\widetilde{\sigma}$ is well-defined. When $d\in 2\N$, not only by hypothesis the sequence $\{s_{dk}\}_{k\in\N_0}$ is positive definite, but also the sequence $\{s_{dk}\}_{k\in\N}$ is positive definite. By Theorem \ref{Prelim:stieltjes} it follows that the points of increase of $\sigma(x)$ are all positive, i.e., $\lambda_i > 0$ for all $i$. This reduces the Hamburger moment sequence to the Stieltjes moment sequence. Hence, $\phi(\lambda_i)$ is well-defined and \eqref{philambdas} holds.
Define the sequence $\{\tilde{s}_k\}_{k\in\N_0}$ by
\begin{equation}
\tilde{s}_k:=\int_{-\infty}^{\infty} x^k d\widetilde{\sigma}(x) \quad \textit{for all }
k\in\N_0.
\end{equation}
We claim that $\{\tilde{s}_k\}_{k\in\N_0}$ is a positive definite completion of the partial sequence $\s$.
Indeed, by the definition of Stieltjes integral one gets that
\begin{equation}\label{series:Stie}
\tilde{s}_k=\sum_{i=0}^{\infty} \mu_i (\phi(\lambda_i))^k.
\end{equation}
Note that the series \eqref{series:Stie} are convergent for each $k\in \N_0$ since
$$
\begin{array}{ll}
\Bigg|\sum_{i=0}^{\infty} \mu_i (\phi(\lambda_i))^k\Bigg|
&\leq \Big|\sum_{|\lambda_i|\geq 1} \mu_i (\phi(\lambda_i))^k\Big| + \Big|\sum_{|\lambda_i| < 1} \mu_i (\phi(\lambda_i))^k\Big| \\
& \leq \Big|\sum_{|\lambda_i|\geq 1} \mu_i \lambda_i^k\Big| + \Big|\sum_{0 < |\lambda_i| < 1} \mu_i \Big| \\
& \leq |s_{kd}|+ s_{0} \\
\end{array}
$$
By (\ref{eq:seq_def_stieltjes2}) it follows that $\tilde{s}_{dk+l_0}=s_{dk+l_0}$ for all $k\in \N_0$.
Hence, $\s$ has a positive definite completion.
\end{proof}
\end{cor}

\begin{remark}
(i) There are positive definite completable patterns that are not of the form as in Corollary \ref{cor:mainresult2}. For example, any pattern $P\subset 2\N_0+1$ is positive definite completable.\\
(ii) If the pattern $P = d\N+l_0$ for $d\in \N$ and $l_0 \in 2\N$, then the pattern $P$ is not positive definite completable. Let $\s$ be a partial positive definite sequence with the pattern $P$. Consider the $2\times2$ principal submatrix
of the partial Hankel matrix $H_{l_0+dm}$ corresponding to $s$ for some $k\in \N_0$
\begin{equation*}
A_m:=
\begin{bmatrix}
s_0 & s_{l_0+dm} \\
s_{l_0+dm} & s_{2l_0+2dm} \\
\end{bmatrix},
\end{equation*}
where $s_0$ is missing and the remaining entries are given. Note that $\det{A_m}>0$ is necessary condition for the pattern $P$ to be positive definite completable. Thus,
\begin{equation*}
s_0>\frac{(s_{l_0+dm})^2}{s_{2(l_0+dm)}}\quad \textit{for all }m\in \N_0.
\end{equation*}
\end{remark}
If the value on the right hand side is increasing with respect to $m$, then there exists no $s_0\in \R$.

\begin{thm}
A pattern $P$ is positive semidefinite completable if and only if
\begin{equation*}
dP +l_0=\{dk+l_0| k\in P\}
\end{equation*}
is positive semidefinite completable for all $d\in \N$ and $l_0 \in 2 \N_{0}$.
\begin{proof}
Let $\s$ be a partial positive semidefinite sequence with the pattern $dP+l_0$.
Then, $\{s_{kd+l_0}\}_{k\in \N_0}$ is a partial positive semidefinite sequence with the pattern $P$.
Since the pattern $P$ is positive semidefinite completable, $\{s_{kd+l_0}\}_{k\in \N_0}$
has a positive semidefinite completion. Then
using the completion $\s$ is a partial positive semidefinite sequence with the pattern $d\N_0+l_0$. By Theorem \ref{thm:mainresult}, $\s$ has a positive semidefinite completion.
Show the converse.
Let $\s$ be a partial positive semidefinite sequence with the pattern $P$. Let $\ts$ be a partial sequence such that $t_{dk+l_0}=s_k$ for all $k\in P$ and the remaining are all missing.
Since $\ts$ is a partial positive semidefinite sequence with the pattern $dP+l_0$, it has a positive semidefinite completion, say $\{\tilde{t}\}_{k\in\N_0}$. Then,
the subsequence $\{\tilde{t}_{dk+l_0}\}_{k\in \N_0}$ is a positive semidefinite completion of $\s$.
\end{proof}
\end{thm}

\begin{cor}
A pattern $P$ is positive definite completable if and only if a pattern $dP$
is positive definite completable for all $d\in \N$.
\end{cor}

By the preceding Theorem and Corollary and by Corollary \ref{cor:trunc} and Remark \ref{remark:semicompletion},
one can check if truncated arithmetic progression patterns are positive completable as follows:

\begin{thm}\label{thm:mainresult3}
Let $P=\{0,1,\ldots,m\}$ and
$d\in \N$, and $l_0 \in 2\N_{0}$.
Then the pattern $dP+l_0$\ is not positive semidefinite completable and
the pattern $dP$ is positive definite completable.
\end{thm}

\section{Positive semidefinite completion}

While quite natural, the connection between positive semidefinite completable and positive definite completable Hankel matrix patterns has not been clarified.  Here we show that if a pattern $ P $ is positive semidefinite completable, then $P $ is positive definite completable. We also characterize certain patterns that do not have positive semidefinite completions, while having positive definite completions.

The following elementary observation is important to characterize positive semidefinite sequences.
\begin{prop}
Suppose $\s$ is a positive semidefinite sequence. Let
$H_n$ be a Hankel matrix corresponding to $\s$.
Let $S:=\{n\in \N~ |~ H_n \textup{ is  singular}\}$.
Then, $S=\{N,N+1,N+2,\cdots\}$ for some $N\in \N_0$. In other words, if $D_0>0$, $D_1>0$, $\cdots$, $D_{N-1}>0$ and $D_{N}=0$ for some $N\in \N$, then $D_{n}=0$ for all $n \geq N$. Here $D_n:=\det{H_n}$.
\begin{proof}
From Definition \ref{defn:definite}, it is clear that $S\neq \emptyset$.
Put $N=\min{\{k:k\in S \}}$. Then $H_n >0$ for all $n < N$ and $H_N \geq 0$, implying that all eigenvalues of $H_{N-1}$ are positive, say that
\begin{equation*}
0 < \lambda_0(H_{N-1}) \leq \lambda_1(H_{N-1}) \leq \lambda_2(H_{N-1}) \leq
\cdots  \leq \lambda_{N-1}(H_{N-1}).
\end{equation*}
Partition $H_N$ as
\begin{equation*}
H_{N}=
\begin{bmatrix}
H_{N-1} & v_{N-1} \\
v_{N-1}^T & s_{2N}
\end{bmatrix},
\end{equation*}
where
$v_n^T=
\begin{bmatrix}
s_{N} & s_{N+1} & \cdots & s_{2N-1}
\end{bmatrix}$.
By Cauchy's Interlace Theorem, one gets that
\begin{equation*}
0 \leq \lambda_0(H_N) \leq \lambda_0(H_{N-1}) \leq \lambda_1(H_N) \leq
\cdots \leq \lambda_{N-1}(H_{N-1}) \leq \lambda_N(H_{N}).
\end{equation*}
Since $H_N\geq 0$, $\lambda_0(H_N)=0$, which means $H_N$ has only one zero eigenvalue.
Moreover, it follows that
\begin{equation*}
0 \leq \cdots \leq \lambda_0 (H_{N+3}) \leq \lambda_0 (H_{N+2}) \leq \lambda_0 (H_{N+1}) \leq \lambda_0 (H_{N})=0,
\end{equation*}
implying $\lambda_0 (H_{n})=0$ for all $n \geq N$.
Therefore, $H_n$ is positive semidefinite and singular for all $n \geq N$.
\end{proof}
\end{prop}

The following is about the positive semidefinite completion of the truncated geometric sequence.
\begin{lem}
Let $a$, $ar$, $ar^2$, $\ldots$ , $ar^{2n}$ be given real numbers with
$n\geq 1$, $a>0$, and $r\in \R$. Then, this sequence has a positive semidefinite completion.
Furthermore, there uniquely exists the completion which is
$s_k=ar^k$ for all $k\in\N_0$.
In fact, if $\sigma(x)$
$(-\infty<x<\infty)$ is a step function with only one point of increase value
$a$ at $x=r$, then
\begin{equation*}
ar^k=\int_{-\infty}^{\infty} x^k d\sigma(x) \quad \textit{for all } k \in\N_0.
\end{equation*}
\begin{proof}
Let $H=H(a,~ar,~\ldots~,~ar^{2n},~?,~?,~?,~?)$ be a partial Hankel matrix with missing entries $s_{2n+1},~s_{2n+2},~s_{2n+3},~s_{2n+4}$, denoted by $?$.
Note the $4\times4$ principal submatrix of $H$
\begin{equation*}
H[\{n,n+1,n+2,n+3\}]=
\begin{bmatrix}
ar^{2n-2} & ar^{2n-1} & ar^{2n} & s_{2n+1} \\
ar^{2n-1} & ar^{2n} & s_{2n+1} & s_{2n+2} \\
ar^{2n} & s_{2n+1} & s_{2n+2} & s_{2n+3}\\
s_{2n+1} & s_{2n+2} & s_{2n+3} & s_{2n+4}
\end{bmatrix}\geq 0
\end{equation*}
if and only if every principal minor of $H$ is nonnegative.
Then it follows that
\begin{equation*}
-(s_{2n+1}-ar^{2n+1})^2 \geq 0,
\end{equation*}
implying that $s_{2n+1}=ar^{2n+1}$. Then, one gets
\begin{equation*}
-r^{2n-2}(s_{2n+2}-ar^{2n+2})^2 \geq 0,
\end{equation*}
so $s_{2n+2}=ar^{2n+2}$.
By induction, the positive semidefinite completion of the given sequence is
$\{ar^k\}_{k\in \N_0}$.
\end{proof}
\end{lem}

\begin{thm}\label{thm:psdpd}
Let $P$ be a pattern of a partial Hankel matrix. If $P$ is positive semidefinite completable, then $P$ is positive definite completable.
\begin{proof}
Suppose that a pattern $P$ is positive semidefinite  completable. Let $A$ be a $(n+1) \times (n+1)$ partial positive definite Hankel matrix with the pattern $P$.
Define the partial Hankel matrix $B(t_0,~t_1,~\ldots~, t_{2n})$ by
\begin{equation*}
t_{k}=
\left\{
\begin{array}{rl}
1/(k+1) & \text{if } k\in P,\\
\textit{missing} & \text{if } k\notin P
\end{array}\right.
\quad\textit{for all } 0\leq k\leq 2n.
\end{equation*}
Choose $\varepsilon>0$ such that $A - \varepsilon B$ is a partial positive definite Hankel matrix.
Since the partial Hankel matrix $A - \varepsilon B$ is defined on the pattern $P$, it has a positive semidefinite completion, say C.
Let $D$ be an $(n+1) \times (n+1)$ Hilbert matrix such that $D_{ij}=1/(i+j+1)$.
Then the $(n+1) \times (n+1)$ Hankel matrix $C+\varepsilon D$ is a positive definite completion of $A$.
\end{proof}
\end{thm}
Note that if a pattern $P$ is not positive definite completable, then $P$ is not positive semidefinite completable.
For example, by Theorem \ref{thm:3x3main} the pattern $P=\{0,1,3,4\}$ of $H_2$ is not positive definite completable. Thus, it is not positive semidefinite completable.

However, the converse of Theorem \ref{thm:psdpd} is not true. By Corollary \ref{cor:trunc} and Remark \ref{remark:semicompletion} the truncated pattern $P$ is positive definite completable, but it is not positive semidefinite completable. Also, the following is another example.

\begin{lem}\label{thm:non-psd}
Suppose that a partial Hankel matrix
\begin{equation*}
\begin{bmatrix}
s_0 &  s_1   & s_2\\
s_1   &  s_2 & s_3\\
s_2 &  s_3 & ?\\
\end{bmatrix}.
\end{equation*}
is partial positive semidefinite. Then,
it has no positive semidefinite completion.
\begin{proof}
Consider the following positive partial Hankel matrix
\begin{equation*}
\begin{bmatrix}
1 &  1   & 1\\
1   &  1 & 2\\
1 &  2 & ?\\
\end{bmatrix}.
\end{equation*}
Then it is partial positive semidefinite, but there is no positive completion.
\end{proof}
\end{lem}
That is, the pattern $P=\{0,1,2,3\}$ of $H_2$ is not positive semidefinite completable. However, by Lemma \ref{lem:oddgiven} the pattern $P=\{0,1,2,3\}$ of $H_2$ is positive definite completable.
In similar, the $P={1,2,3,4}$ of $H_2$ is not positive semidefinite completable.

\section{Positive Hankel matrix completion}

In this section we provide a complete characterization of all $ 3 \times 3 $ partial positive definite Hankel matrices that have positive definite completions.  There are only three patterns that do not admit such completions. For $ 4 \times 4 $ partial positive Hankel matrices, there are many more patterns that are not completable. While over 30 of these patterns have been checked a complete list of non-completable patterns (among 128 possible patterns) does not follow any obvious symmetry. Using these two cases, the last theorem of this section extends these sets of results to all $ n \times n $ matrices.

We begin with $3 \times 3$ Hankel matrix completions.

\begin{lem}\label{lem:3x3ex1}
Suppose that a partial Hankel matrix
\begin{equation*}
H=
\begin{bmatrix}
s_0 &  s_1   & s_2\\
s_1   &  s_2 & ?\\
s_2 &  ? & s_4\\
\end{bmatrix}.
\end{equation*}
is partial positive definite. Then,
it has a positive definite completion.
Furthermore, the missing entry $s_3$ satisfies the following inequality:
\begin{equation*}
\frac{-\sqrt{P_1P_2}+s_1 s_2}{s_0} < s_3 < \frac{\sqrt{P_1P_2}+s_1 s_2}{s_0},
\end{equation*}
where
\begin{align*}
P_1=\det{
\begin{bmatrix}
s_0 &  s_1 \\
s_1   &  s_2 \\
\end{bmatrix}}
\textit{ and  }P_2=\det{
\begin{bmatrix}
s_0 &  s_2 \\
s_4   &  s_4 \\
\end{bmatrix}}.
\end{align*}
In addition, if partial Hankel matrices
\begin{equation*}
\begin{bmatrix}
s_0 &  ?  & s_2\\
?   &  s_2 & ?\\
s_2 &  ? & s_4\\
\end{bmatrix}
\text{ and }
\begin{bmatrix}
s_0 &  ?  & s_2\\
?   &  s_2 & ?\\
s_2 &  ? & ?\\
\end{bmatrix}.
\end{equation*}
are partial positive definite, then
they have positive definite completions.
\begin{proof}
Since the Hankel matrix $H$ is partial positive definite,
$P_1 > 0$ and $P_2 >0$.
The inequality comes from $\det{H_2}$. Note that $H_2>0$ if and only if all leading principal minors are positive.
\end{proof}
\end{lem}

\begin{remark}\label{remark:sylvester}
The following fact is useful:
\begin{equation*}
\begin{bmatrix}
s_0 &  s_1   & s_2\\
s_1   &  s_2 & s_3\\
s_2 &  s_3 & s_4\\
\end{bmatrix}>0
\textup{ iff }
\begin{bmatrix}
s_4 &  s_3   & s_2\\
s_3   &  s_2 & s_1\\
s_2 &  s_1 & s_0\\
\end{bmatrix}>0,
\end{equation*}
since two matrices are unitarily similar. For example, the following partial Hankel matrices can be checked  whether they have positive definite completions in the same way.
\begin{equation*}\label{3by3:ex1}
\begin{bmatrix}
s_0 &  ?   & s_2\\
?   &  s_2 & s_3\\
s_2 &  s_3 & s_4\\
\end{bmatrix} \textup{ and }
\begin{bmatrix}
s_4 & s_3   & s_2\\
s_3   &  s_2 & ?\\
s_2 & ? & s_0\\
\end{bmatrix}
\end{equation*}
\end{remark}

\begin{thm}\label{thm:3x3main}
For $3\times3$ partial Hankel matrices $H_2$,
the following patterns are not positive definite completable.
\begin{equation*}
\{0,1,3,4\},~ \{0,1,4\},~\{0,3,4\}
\end{equation*}
\begin{proof}
Consider the following positive partial Hankel matrix
\begin{equation*}\label{3by3:counterex1}
H=
\begin{bmatrix}
1 &  \frac{1}{2}   & ?\\
\frac{1}{2}   &  ? & ?\\
? &  ? & \frac{1}{16}\\
\end{bmatrix}.
\end{equation*}
Since the $2\times 2$ principal submatrices $H[\{1,2\}]$ and $H[\{1,3\}]$ cannot have positive definite completions at the same time, $H$ does not have positive definite completion.
\end{proof}
\end{thm}

\begin{remark}
In fact, all patterns $P\subset\{0,1,2,3,4\}$ of a partial Hankel matrix $H_2$ are positive definite completable except the above three patterns.
\end{remark}
Now consider $4 \times 4$ Hankel matrix completions.

\begin{thm}\label{thm:4x4main}
For $4\times4$ partial Hankel matrices
the following patterns is not positive definite completable:
\begin{align*}
P=& \{0,1,2,4,5,6\}
\end{align*}
\end{thm}


\begin{thm}
If a partial positive definite Hankel matrix has a $3\times3$ or $4\times4$ principal submatrix with one of the previous pattern which is not positive definite Hankel completable, then the matrix has no positive definite Hankel completion.
\end{thm}

\begin{ex}
The partial positive Hankel matrix
\begin{equation*}
H=\begin{bmatrix}
s_0 &  s_1 & ? & s_3 & s_4\\
s_1 &  ? & s_3 & s_4 & ?\\
? &  s_3 & s_4 & ? & s_6\\
s_3 &  s_4 & ? & s_6 & s_7\\
s_4 &  ? & s_6 & s_7 & ?\\
\end{bmatrix}
\end{equation*}
may not have a positive definite completion since by Theorem $\ref{thm:3x3main}$ the $3\times3$ principal submatrix $H[{1,2,3}]$ may not have a positive definite completion.
\end{ex}

\section*{Acknowledgment}
The authors wish to express their gratitude to the
anonymous referees for their careful reading of the manuscript and their helpful suggestions.


\end{document}